%% file: CSFEM_Flut.tex
\documentclass[onecolumn]{revtex4}
\usepackage{latexsym}
\usepackage{amsmath,amsthm}
\usepackage{times}
\usepackage{amssymb}
\usepackage{fancyheadings}
\usepackage[T1]{fontenc}
\usepackage[utf8]{inputenc}
\usepackage{fancyhdr}
\usepackage{url}
\usepackage{hyperref}
\usepackage{color}
\usepackage{subfigure}
\usepackage{multirow}
\usepackage{graphicx}

\newcommand{\xx}{\mathbf{x}}
\newcommand{\KK}{\mathbf{K}}
\newcommand{\bfm}{\mathbf{M}}
\newcommand{\bn}{\mathbf{N}}
\newcommand{\DD}{\mathbf{D_b}}

\newcommand{\bveps}{\boldsymbol{\varepsilon}}
\newcommand{\uu}{\mathbf{u}}
\newcommand{\BB}{\mathbf{B}}

\newcommand{\pp}{\mathbf{p}}
\newcommand{\bss}{\mathbf{s}}
\newcommand{\bdb}{\mathbf{b}}

\newcommand{\bds}{\boldsymbol{\delta}}

\newcommand{\rmd}{\mathrm{d}}

\newcommand\undermat[2]{%
  \makebox[0pt][l]{$\smash{\underbrace{\phantom{%
    \begin{matrix}#2\end{matrix}}}_{\text{$#1$}}}$}#2}

\theoremstyle{remark}

\newcommand{\Eref}[1]{Equation (\ref{#1})}
\newcommand{\fref}[1]{Figure (\ref{#1})}

\begin{document}


\title{\bf Linear flutter analysis of functionally graded panels using cell based smoothed finite element method and discrete shear gap technique}
\author{Sundararajan Natarajan$^{1}$\vspace*{3mm}}\thanks{Corresponding author. Email: sundararajan.natarajan@gmail.com,  Tel.: +61 293855030.}
\author{Karthik Kaleeswaran$^2$\vspace*{3mm}}
\author{Ganapathi Manickam$^3$\vspace*{3mm}}
\affiliation{$^1$School of Civil and Environmental Engineering,
University of New South Wales, Sydney, Australia \vspace*{3mm}}
\affiliation{$^2$Siemens India Ltd., Bangalore, India \vspace*{3mm}}
\affiliation{$^3$Stress \& DTA, IES-Aerospace, Mahindra Satyam Computer Services Ltd., Bangalore, India \vspace*{3mm}}
\date{\today}


\begin{abstract}
\noindent In this paper, a cell-based smoothed finite element method with discrete shear gap technique for triangular elements is employed to study the linear flutter characteristics of functionally graded material (FGM) flat panels. The influence of thermal environment, the presence of a centrally located circular cutout and the aerodynamic damping on the supersonic flutter characteristics of flat FGM panels is also investigated. The structural formulation is based on the first-order shear deformation theory and the material properties are assumed to be temperature dependent and graded only in the thickness direction according to power law distribution in terms of the volume fraction of its constituent materials. The aerodynamic force is evaluated by considering the first order high mach number approximation to linear potential flow theory. The formulation includes transverse shear deformation and in-plane and rotary inertia effects. The influence of the plate thickness, aspect ratio, boundary conditions, material gradient index, temperature dependent material properties, damping, cutout size, skewness of the plate and boundary conditions on the critical aerodynamic pressure is numerically studied.

\vspace*{2ex}\noindent\textit{\bf Keywords}: Cell-based smoothed finite element method, Discrete shear gap technique, Functionally graded material, Material gradient index, Flutter.
\end{abstract}

\maketitle

\thispagestyle{fancy}

\section{Introduction}
\label{Intro}

\vspace*{0.2cm}
In recent years, a new class of engineered material, the functionally graded materials (FGMs) has attracted researchers to investigate its structural behaviour. The FGMs are microscopically inhomogeneous and the mechanical and the thermal properties vary \textit{smoothly and continuously} from one surface to another. FGMs combine the best properties of its constituents. Typically, the FGMs are made from a mixture of ceramic and metal. The ceramic constituent provides thermal stability due to its low thermal conductivity, whilst the metallic constituent provides structural stability. FGMs eliminate the sharp interfaces existing in laminated composites with a gradient interface and are considered to be an alternative in many engineering sectors such as the aerospace industry, biomechanics industry, nuclear industry, tribology, optoelectronics and other high performance applications where the structural member is exposed to high thermal gradient in addition to mechanical load. 

The static and the dynamic characteristics have been studied in detail by many researchers using different plate theories, for example, first order shear deformation theory (FSDT)~\cite{Reddy2000,Yang2002,Sundararajan2005}, second and other higher order accurate theory~\cite{Qian2004a,Ferreira2006,natarajanmanickam2012} have been used to describe the plate kinematics. Existing approaches in the literature to study plate and shell structures made up of FGMs uses finite element method (FEM) based on Lagrange basis functions~\cite{Reddy2000,ganapathiprakash2006,Sundararajan2005}, meshfree methods~\cite{Qian2004a,Ferreira2006} and recently Valizadeh \textit{et al.,}~\cite{valizadehnatarajan2013} used non-uniform rational B-splines based FEM to study the static and the dynamic characteristics of FGM plates in thermal environment. Akbari \textit{et al.,}~\cite{R2010} studied two-dimensional wave propagation in functionally graded solids using the meshless local Petrov-Galerkin method. Huang \textit{et al.,}~\cite{Huang2011} proposed solutions for the free vibration of side-cracked FGM thick plates based on Reddy's third-order shear deformation theory using Ritz technique. Kitipornchai~\textit{et al.,}~\cite{Kitipornchai2009} studied nonlinear vibration of edge cracked functionally graded Timoshenko beams using Ritz method. Yang~\textit{et al.,}~\cite{Yang2010} studied the nonlinear dynamic response of a functionally graded plate with a through-width crack based on Reddy's third-order shear deformation theory using a Galerkin method. Natarajan \textit{et al.,}~\cite{natarajanbaiz2011,natarajanbaiz2011a} and Baiz \textit{et al.,}~\cite{baiznatarajan2011} studied the influence of the crack length on the free flexural vibrations and on the critical buckling load of FGM plates using the XFEM and smoothed XFEM, respectively. Plates with cutouts are extensively used in transport vehicle structures. Cutouts are made to lighten the structure, for ventilation, to provide accessibility to other parts of the structures and for altering the resonant frequency. Therefore, the natural frequencies of plates with cutouts are of considerable interest to designers of such structures. Most of the earlier investigations on plates with cutouts have been confined to isotropic plates~\cite{paramasivam1973,aliatwal1980,huangsakiyama1999} and laminated composites~\cite{reddy1982,sivakumariyengar1998}. Recently, Janghorban and Zare~\cite{janghorbanzare2011} studied the influence of cutout on the fundamental frequency of FGM plates in thermal environment using FEM. Recently Akbari \textit{et al.,}~\cite{rahimabadinatarajan2013} employed the XFEM to study the influence of internal discontinuities, viz., cracks and cutouts on the fundamental frequency of FGM plates in thermal environment. The presence of a cutout can also influence the flutter characteristics. Futhermore, in practice, the use of these materials in aerospace industries has necessitated to understand the dynamic characteristics of functionally graded structures. This has attracted researchers~\cite{prakashganapathi2006,navazihaddadpour2007,ibrahimtawfik2007} to study the flutter characteristics of FGM panels. The above list is no way comprehensive and interested readers are referred to the literature and references therein and a recent review paper by Jha and Kant~\cite{jhakant2013} on FGM plates. 

\vspace*{0.2cm}
\noindent
{\bf Approach} 
In this paper, we study the linear flutter characteristics of FGM flat panels using a 3-noded triangular element. A cell-based smoothing technique combined with discrete shear gap method is employed for this study. The influence of the plate thickness, aspect ratio, boundary conditions, material gradient index, temperature dependent material properties, damping, cutout size, skewness of the plate and boundary conditions on the critical aerodynamic pressure is numerically studied. 

\vspace*{0.2cm}
\noindent
{\bf Approach} The paper is organized as follows, the next section will give an introduction to FGM and a brief overview of Reissner-Mindlin plate theory. Section~\ref{csdsg3} describes the cell-based smoothing technique combined with discrete shear gap method for 3-noded triangular elements. The efficiency of the present formulation, numerical results and parametric studies are presented in Section~\ref{numres}, followed by concluding remarks in the last section.

\vspace*{0.2cm}
\section{Theoretical Formulation}\label{theory}

\input{fgm}

\input{plateTheory}

\vspace*{0.2cm}
\section{Spatial discretization}\label{csdsg3}
In this study, three-noded triangular element with five degrees of freedom (dofs) $\boldsymbol{\delta} = \{u,v,w,\theta_x,\theta_y\}$ is employed. The displacement field is approximated by
\begin{equation}
\uu^h = \sum_IN_I \boldsymbol{\delta}_I
\end{equation}
where $\boldsymbol{\delta}_I$ are the nodal dofs and $N_I$ are the standard finite element shape functions given by
\begin{equation}
N = \left[ 1-\xi-\eta, \;\; \eta, \;\; \xi \right]
\end{equation}

\begin{figure}[htpb]
\centering
\scalebox{0.8}{\input{./Figures/triangle.pstex_t}}
\caption{A triangular element is divided into three subtriangles. $\Delta_1, \Delta_2$ and $\Delta_3$ are the subtriangles created by connecting the central point $O$ with three field nodes.}
\label{fig:triEle}
\end{figure}

In the proposed approach, cell-based smoothed finite element method (CSFEM) is combined with stabilized discrete shear gap method (DSG) for three-noded triangular element, called as `cell-based discrete shear gap method (CS-DSG3).' The cell-based smoothing technique decreases the computational complexity, whilst DSG suppresses the shear locking phenomenon when the present formulation is applied to thin plates. Interested readers are referred to the literature and references therein for the description of cell-based smoothing technique~\cite{liudai2007,bordasnatarajan2010} and DSG method~\cite{bletzingerbischoff2000}. In the CS-DSG3, each triangular element is divided into three subtriangles. The displacement vector at the center node is assumed to be the simple average of the three displacement vectors of the three field nodes. In each subtriangle, the stabilized discrete shear gap (DSG3) (Note: 3 refers to discrete shear gap technique applied to 3-noded triangular element) is used to compute the strains and also to avoid the transverse shear locking. Then the strain smoothing technique on the whole triangular element is used to smooth the strains on the three subtriangles. Consider a typical triangular element $\Omega_e$ as shown in \fref{fig:triEle}. This is first divided into three subtriangles $\Delta_1, \Delta_2$ and $\Delta_3$ such that $\Omega_e = \bigcup\limits_{i=1}^3 \Delta_i$. The coordinates of the center point $\xx_o = (x_o,y_o)$ is given by:
\begin{equation}
(x_o,y_o) = \frac{1}{3} (x_I, y_I)
\end{equation}
The displacement vector of the center point is assumed to be a simple average of the nodal displacements as
\begin{equation}
\boldsymbol{\delta}_{eO} = \frac{1}{3}\boldsymbol{\delta}_{eI}
\label{eqn:centerdefl}
\end{equation}
The constant membrane strains, the bending strains and the shear strains for subtriangle $\Delta_1$ is given by:
\begin{align}
\bveps_p &= \left[ \begin{array}{ccc} \pp_1^{\Delta_1} & \pp_2^{\Delta_1} & \pp_3^{\Delta_1} \end{array} \right] \left\{ \begin{array}{c} \bds_{eO} \\ \bds_{e1} \\ \bds_{e2} \end{array} \right\} \nonumber \\
\bveps_b &= \left[ \begin{array}{ccc} \bdb_1^{\Delta_1} & \bdb_2^{\Delta_1} & \bdb_3^{\Delta_1} \end{array} \right] \left\{ \begin{array}{c} \bds_{eO} \\ \bds_{e1} \\ \bds_{e2} \end{array} \right\} \nonumber \\
\bveps_s &= \left[ \begin{array}{ccc} \bss_1^{\Delta_1} & \bss_2^{\Delta_1} & \bss_3^{\Delta_1} \end{array} \right] \left\{ \begin{array}{c} \bds_{eO} \\ \bds_{e1} \\ \bds_{e2} \end{array} \right\}
\label{eqn:constrains}
\end{align}
Upon substituting the expression for $\boldsymbol{\delta}_{eO}$ in \Eref{eqn:constrains}, we obtain:
\begin{align}
\bveps_p ^{\Delta_1}&= \left[ \begin{array}{ccc} \frac{1}{3}\pp_1^{\Delta_1} + \pp_2^{\Delta_1}& \frac{1}{3}\pp_1^{\Delta_1} + \pp_3^{\Delta_1} & \frac{1}{3} \pp_1^{\Delta_1} \end{array} \right] \left\{ \begin{array}{c} \bds_{e1} \\ \bds_{e2} \\ \bds_{e3} \end{array} \right\} = \BB_p^{\Delta_1} \bds_e \nonumber \\
\bveps_b^{\Delta_1} &= \left[ \begin{array}{ccc} \frac{1}{3}\bdb_1^{\Delta_1} + \bdb_2^{\Delta_1}& \frac{1}{3}\bdb_1^{\Delta_1} + \bdb_3^{\Delta_1} & \frac{1}{3} \bdb_1^{\Delta_1} \end{array} \right] \left\{ \begin{array}{c} \bds_{e1} \\ \bds_{e2} \\ \bds_{e3} \end{array} \right\} = \BB_b^{\Delta_1} \bds_e\nonumber \\
\bveps_s^{\Delta_1} &= \left[ \begin{array}{ccc} \frac{1}{3}\bss_1^{\Delta_1} + \bss_2^{\Delta_1}& \frac{1}{3}\bss_1^{\Delta_1} + \bss_3^{\Delta_1} & \frac{1}{3} \bss_1^{\Delta_1} \end{array} \right] \left\{ \begin{array}{c} \bds_{e1} \\ \bds_{e2} \\ \bds_{e3} \end{array} \right\} = \BB_s^{\Delta_1} \bds_e \nonumber \\
\label{eqn:constrainsSubtri}
\end{align}
where $\pp_i, (i=1,2,3)$, $\bdb_i, (i=1,2,3)$ and $\bss_i, (i=1,2,3)$ are given by:
\begin{align}
\BB_p &= \frac{1}{2A_e} \left[ \begin{array}{rrrrrrrrrrrrrrr} b-c & 0 & 0 & 0 & 0 & c & 0 & 0 & 0 & 0 & -b & 0 & 0 & 0 & 0 \\ 0 & d-a & 0 & 0 & 0 & 0 & -d & 0 & 0 & 0 & a & 0 & 0 & 0 & 0 \\\undermat{\pp_1}{ d-a & b-c & 0 & 0 & 0} & \undermat{\pp_2}{-d & c & 0 & 0 & 0} & \undermat{\pp_3}{a & -b & 0 & 0 & 0} \end{array} \right] \nonumber \\ \nonumber \\ \nonumber \\
\BB_b & = \frac{1}{2A_e} \left[ \begin{array}{rrrrrrrrrrrrrrr} 0 & 0 & 0 & b-c & 0 & 0 & 0 & 0 & c & 0 & 0 & 0 & 0 & -b & 0 \\ 0 & 0 & 0 & 0 & d-a & 0 & 0 & 0 & 0 & -d & 0 & 0 & 0 & 0 & a \\ \undermat{\bdb_1}{0 & 0 & 0 & d-a & b-c} & \undermat{\bdb_2}{0 & 0 & 0 & -d & c} & \undermat{\bdb_3}{0 & 0 & 0 & a & -b} \end{array} \right] \nonumber \\ \nonumber \\ \nonumber \\
\BB_s & = \frac{1}{2A_e} \left[ \begin{array}{rrrrrrrrrrrrrrr} 0 & 0 & b-c & A_e & 0 & 0 & 0 & c & ac/2 & bc/2 & 0 & 0 & -b & -bd/2 & -bc/2 \\ \undermat{\bss_1}{0 & 0 & d-a & 0 & A_e} & \undermat{\bss_2}{0 & 0 & -d & -ad/2 & -bd/2} & \undermat{\bss_3}{0 & 0 & a & ad/2 & ac/2} \end{array} \right] \\
\end{align}
where $a = x_2 - x_1; b = y_2 - y_1; c = y_3 - y_1$ and $d = x_3 - x_1$ (see \fref{fig:dsg3}), $A_e$ is the area of the triangular element and $\BB_s$ is altered shear strains. The strain-displacement matrix for the other two triangles can be obtained by cyclic permutation.
\begin{figure}[htpb]
\centering
\scalebox{0.7}{\input{./Figures/stdtri.pstex_t}}
\caption{Three-noded triangular element and local coordinates in discrete shear gap method.}
\label{fig:dsg3}
\end{figure}
Now applying the cell-based strain smoothing~\cite{liudai2007,bordasnatarajan2010}, the constant membrane strains, the bending strains and the shear strains are respectively employed to create a smoothed membrane strain $\overline{\bveps}_p$, smoothed bending strain $\overline{\bveps}_b$ and smoothed shear strain $\overline{\bveps}_s$on the triangular element $\Omega_e$ as:

\begin{align}
\overline{\bveps}_p &= \int\limits_{\Omega_e} \bveps_b \Phi_e(\xx)~\rmd \Omega = \sum\limits_{i=1}^3 \bveps_p^{\Delta_i} \int\limits_{\Delta_i} \Phi_e(\xx)~\rmd \Omega \nonumber \\
\overline{\bveps}_b &= \int\limits_{\Omega_e} \bveps_b \Phi_e(\xx)~\rmd \Omega = \sum\limits_{i=1}^3 \bveps_b^{\Delta_i} \int\limits_{\Delta_i} \Phi_e(\xx)~\rmd \Omega \nonumber \\
\overline{\bveps}_s &= \int\limits_{\Omega_e} \bveps_s \Phi_e(\xx)~\rmd \Omega = \sum\limits_{i=1}^3 \bveps_s^{\Delta_i} \int\limits_{\Delta_i} \Phi_e(\xx)~\rmd \Omega
\end{align}
where $\Phi_e(\xx)$ is a given smoothing function that satisfies. In this study, following constant smoothing function is used:
\begin{equation}
\Phi(\xx) = \left\{ \begin{array}{cc} 1/A_c & \xx \in \Omega_c \\ 0 & \xx \notin \Omega_c \end{array} \right.
\end{equation}
where $A_c$ s the area of the triangular element, the smoothed membrane strain, the smoothed bending strain and the smoothed shear strain is then given by
\begin{equation}
\left\{ \overline{\bveps}_p, \overline{\bveps}_b, \overline{\bveps}_s \right\} = \frac{ \sum\limits_{i=1}^3 A_{\Delta_i} \left\{ \bveps_p^{\Delta_i}, \bveps_b^{\Delta_i}, \bveps_s^{\Delta_i} \right \} }{A_e}
\end{equation}
The smoothed elemental stiffness matrix is given by
\begin{align}
\KK &= \int\limits_{\Omega_e} \overline{\BB}_p \mathbf{A} \overline{\BB}_p^{\rm T} + \overline{\BB}_p \BB \overline{\BB}_b^{\rm T} + \overline{\BB}_b \BB \overline{\BB}_p^{\rm T} + \overline{\BB}_b \mathbf{D} \overline{\BB}_b^{\rm T} + \overline{\BB}_s \mathbf{E} \overline{\BB}_s^{\rm T}~\rmd \Omega \nonumber \\
&= \left( \overline{\BB}_p \mathbf{A} \overline{\BB}_p^{\rm T} + \overline{\BB}_p \BB \overline{\BB}_b^{\rm T} + \overline{\BB}_b \BB \overline{\BB}_p^{\rm T} + \overline{\BB}_b \mathbf{D} \overline{\BB}_b^{\rm T} + \overline{\BB}_s \mathbf{E} \overline{\BB}_s^{\rm T} \right) A_e
\end{align}
where $\overline{\BB}_p, \overline{\BB}_b$ and $\overline{\BB}_s$ are the smoothed strain-displacement matrix. The mass matrix $\bfm$, the geometric stiffness matrix $\KK_G$ and the aerodynamic matrices $\overline{\mathbf{A}}$ and $\mathbf{D}_A$ are computed by following the conventional finite element procedure.

\vspace*{0.2cm}
\section{Numerical Results}\label{numres}
In this section, we present the critical aerodynamic pressure and the critical frequency of functionally graded material plates immersed in a supersonic flow using three-noded triangular element with cell-based smoothed finite element method and discrete shear gap technique. The element has five degrees of freedom $(u_o,v_o,w_o,\theta_x,\theta_y)$. The shear locking phenomenon is suppressed with a combination of the discrete shear gap technique and the strain smoothing method. The FGM plate considered here is made up of silicon nitride (Si$_3$N$_4$) and stainless steel (SUS304). The material is considered to be temperature dependent and the temperature coefficients corresponding to Si$_3$N$_4$/SUS304 are listed in Table \ref{table:tempdepprop}~\cite{Reddy1998,Sundararajan2005}. The mass density $(\rho)$ and the thermal conductivity $(\kappa)$ are $\rho_c=$ 2370 kg/m$^3$, $\kappa_c=$ 9.19 W/mK for Si$_3$N$_4$ and $\rho_m=$ 8166 kg/m$^3$, $\kappa_m=$ 12.04 W/mK for SUS304. Poisson's ratio $\nu$ is assumed to be constant and taken as 0.28 for the current study~\cite{Sundararajan2005}. Here the modified shear correction factor obtained based on energy equivalence principle as outlined in~\cite{Singh2011} is used. The boundary conditions for simply supported and clamped cases are: \\
\emph{Simply supported boundary condition}: 
\begin{equation}
u_o = w_o = \theta_y = 0 \hspace{0.2cm} \textup{on} \hspace{0.2cm} x=0,a; \hspace{0.2cm} v_o = w_o = \theta_x = 0 \hspace{0.2cm} \textup{on} \hspace{0.2cm} y = 0,b
\end{equation}
\emph{Clamped boundary condition}:
\begin{equation}
u_o = v_o = w_o = \theta_x = \theta_y= 0 \hspace{0.2cm} \textup{on} \hspace{0.2cm} x = 0, a \hspace{0.2cm} \& \hspace{0.2cm} y = 0,b
\end{equation}

\begin{table}
\renewcommand\arraystretch{1.5}
\caption{Temperature dependent coefficient for material Si$_3$N$_4$/SUS304, Ref~\cite{Reddy1998,Sundararajan2005}.}
\centering
\begin{tabular}{lcccccc}
\hline
Material & Property & $P_o$ & $P_{-1}$ & $P_1$ & $P_2$ & $P_3$  \\
\hline
\multirow{2}{*}{Si$_3$N$_4$} & $E$(Pa) & 348.43e$^9$ &0.0& -3.070e$^{-4}$ & 2.160e$^{-7}$ & -8.946$e^{-11}$  \\
& $\alpha$ (1/K) & 5.8723e$^{-6}$ & 0.0 & 9.095e$^{-4}$ & 0.0 & 0.0 \\
\cline{2-7}
\multirow{2}{*}{SUS304} & $E$(Pa) & 201.04e$^9$ &0.0& 3.079e$^{-4}$ & -6.534e$^{-7}$ & 0.0  \\
& $\alpha$ (1/K) & 12.330e$^{-6}$ & 0.0 & 8.086e$^{-4}$ & 0.0 & 0.0 \\
\hline
\end{tabular}
\label{table:tempdepprop}
\end{table}

\paragraph*{Skew boundary transformation} For skew plates, the edges of the boundary elements may not be parallel to the global axes $(x,y,z)$. In order to specify the boundary conditions on skew edges, it is necessary to use the edge displacements $(u_o^\prime,v_o^\prime,w_o^\prime)$ etc, in a local coordinate system $(x^\prime,y^\prime,z^\prime)$ (see \fref{fig:platefig}). The element matices corresponding to the skew edges are transformed from global axes to local axes on which the boundary conditions can be conveniently specified. The relation between the global and the local degrees of freedom of a particular node is obtained by:
\begin{equation}
\boldsymbol{\delta} = \mathbf{L}_g \boldsymbol{\delta}^\prime
\end{equation}
where $\boldsymbol{\delta}$ and $\boldsymbol{\delta}^\prime$ are the generalized displacement vector in the global and the local coordinate system, respectively. The nodal transformation matrix for a node $I$ on the skew boundary is given by:
\begin{equation}
\mathbf{L}_g = \left[ \begin{array}{rrrrr} \cos\psi & \sin\psi & 0 & 0 & 0 \\ -\sin\psi & \cos\psi & 0 & 0 & 0 \\ 0 & 0 & 1 & 0 & 0 \\ 0 & 0 & 0 & \cos\psi & \sin\psi \\ 0 & 0 & 0 & -\sin\psi & \cos\psi \end{array} \right]
\end{equation}
where $\psi$ defines the skewness of the plate.

\paragraph*{Validation} Before proceeding with the detailed study, the formulation developed herein is validated against available results pertaining to the critical aerodynamic pressure and the critical frequency for isotropic plate and functionally graded material plates. The computed aerodynamic pressure and frequency: (a) for an isotropic plate immersed in normal flow is given in Table~\ref{fig:meshconvestudy} and (b) for functionally graded material in thermal environment immersed in a normal flow is given in Table~\ref{table:fgmtempdependflut}. Based on a progressive refinement, a 40 $\times$ 40 structured triangular mesh is found to be adequate to model the full plate. The results evaluated for both simply supported and clamped boundary conditions are found to be in very good agreement with the results in the literature~\cite{liaosun1993,chowdarysinha1996,singhaganapathi2005,prakashganapathi2006}.

\begin{table}[htpb]
\renewcommand\arraystretch{1}
\caption{Mesh convergence study of critical aerodynamic pressure $\lambda_{cr} = \lambda a^3/(\pi^4D)$ of isotropic plates $(a/b=$ 1, $a/h=$ 100$)$.}
\centering
\begin{tabular}{lrr}
\hline
Mesh & \multicolumn{2}{r}{Skew angle $\psi$}\\
\cline{2-3}
& 0$^\circ$ & 30$^\circ$ \\
\hline
8 $\times$ 8 & 6.7391 & 7.2042 \\
16 $\times$ 16 & 5.4478 & 6.6268 \\
32 $\times$ 32 & 5.2954 & 6.4984 \\
40 $\times$ 40 & 5.2794 & 6.4824 \\
Ref.~\cite{singhaganapathi2005} & 5.27 & 6.47 \\
Ref.~\cite{liaosun1993} & 5.12 & 6.31 \\
Ref.~\cite{chowdarysinha1996} & 5.25 & 6.82 \\
\hline
\end{tabular}
\label{fig:meshconvestudy}
\end{table}

\begin{table}[htpb]
\renewcommand\arraystretch{1.5}
\caption{Comparison of flutter behaviour of temperature dependent FGM plate (material: Si$_3$N$_4$/SUS304, $a/h=$ 20, all edges simply supported.}
\centering
\begin{tabular}{lcccccc}
\hline
$T_c,T_m$ & gradient index & \multicolumn{2}{c}{$\overline{\omega}_{cr}^{2\ast}$} & & \multicolumn{2}{c}{$\lambda_{cr}$} \\
\cline{3-7}
 & $n$ & Ref.~\cite{prakashganapathi2006} & Present && Ref.~\cite{prakashganapathi2006} & Present \\
 \hline
 \multirow{3}{*}{300,300} & 0 & 9661.35 & 9653.20 && 775.78 & 775.98 \\
 & 1 & 3515.57 & 3474.40 && 625.78 & 618.95 \\
 & 5 & 2348.72 & 2326.20 && 571.48 & 566.60 \\
 \cline{2-7}
 \multirow{3}{*}{600,300} & 0 & 7475.77 & 7470.50 && 647.65 & 647.85 \\
 & 1 & 2528.99 & 2520.10 && 499.61 & 496.29 \\
 & 5 & 1554.78 & 1547.70 && 433.20 & 430.66\\
 \hline
\end{tabular}
\label{table:fgmtempdependflut}
\end{table}

\begin{table}[htpb]
\renewcommand\arraystretch{1.5}
\caption{Influence of the plate aspect ratio and the material gradient index $n$ on the flutter behaviour of Si$_3$N$_4$/SUS304 with $a/h=$ 100. The ceramic property is used for normalization.}
\centering
\begin{tabular}{ccrrrrrr}
\hline
 & $a/b$ & \multicolumn{6}{c}{gradient index $n$}\\
\cline{3-8} 
 &  & 0 & 1 & 2 & 3 & 4 & 5 \\
\hline
\multirow{5}{*}{$\lambda_{cr}$} & 0.5 & 357.0313 & 285.1563 & 273.4375 & 267.9688 & 264.0625 & 260.9375 \\
& 1.0 & 476.5625 & 380.4688 & 364.0625 & 357.0313 & 352.3438 & 348.4375 \\
& 2.0 & 1026.5625 & 817.1875 & 782.0313 & 767.1875	 & 757.0313	 & 749.2188 \\
& 3.0 & 2160.1563 & 1702.3438 & 1628.9063 & 1600.7813 & 1582.0313	 & 1567.1875 \\
& 5.0 & 7170.3125 & 5224.2188 & 4995.3125 & 4957.0313 & 4949.2188	 & 4946.0938 \\
\cline{2-8}
\multirow{5}{*}{$\overline{\omega}_{cr}^2$} & 0.5	& 899.5877 & 323.7708 & 262.5403 & 238.7031 & 225.3648 & 216.6146 \\
& 1.0 & 1714.4392	 & 617.5253	 & 499.4822 & 454.2823 & 	429.4400 &	413.1416 \\
& 2.0 & 5380.1598 & 1939.1894 & 1568.6419 & 1426.4056 & 1347.8699 & 1296.8628 \\
& 3.0 & 16296.3293 & 5848.9675 & 4731.1698 & 4305.1735 & 4070.6429 & 3917.9894 \\
& 5.0 & 91276.9666 & 31458.0844 & 25439.0947 & 23246.8634 & 22084.4147 & 21346.8356 \\
 \hline
\end{tabular}
\label{table:plateaspectratio}
\end{table}

Table \ref{table:plateaspectratio} presents the flutter characteristics of square and rectangular plates with $a/h=$ 100 made up of Si$_3$N$_4$/SUS304 is investigated, neglecting the influence of thermal load. It is inferred from Table \ref{table:plateaspectratio} that the critical aerodynamic pressure decreases with increase in the material gradient index $n$. However, the rate of decrease of flutter speed is high for low value of $n$. This can be attributed to the increase in the metallic volume fraction. Also, it can be observed that, for a given thickness ratio, the critical aerodynamic pressure increases with the increase in the plate aspect ratio $a/b$.

\begin{table}[htpb]
\renewcommand\arraystretch{1.5}
\caption{Flutter behaviour of temperature dependent FGM plate (material: Si$_3$N$_4$/SUS304) with $a/b=$ 1 and $a/h=$ 20.}
\centering
\begin{tabular}{llrrrrrrrr}
\hline
$T_c,T_m$ & gradient index & \multicolumn{2}{c}{In-vacuo}& & \multicolumn{2}{c}{Without aerodynamic damping} & & \multicolumn{2}{c}{With aerodynamic damping}\\
\cline{3-4} \cline{6-7} \cline{9-10}
 & $n$ & $\overline{\omega}_{1}^2$ & $\overline{\omega}_{2}^2$ & & $\overline{\omega}_{cr}^2$ & $\lambda_{cr}$ & & $\overline{\omega}_{cr}^2$ & $\lambda_{cr}$\\
\hline
\multirow{6}{*}{300,300} & 0 & 2051.9860 & 12542.3513 & &  9653.0162 & 776.5625 & & 9746.1925 & 787.8125\\
& 1 & 746.2407 & 4524.0162 & & 3474.4000 & 619.5313 & & 3529.8134 & 635.7813\\
& 2 & 603.7262 & 3658.6725 & &2811.2893 & 592.9688 & & 2870.6620 & 614.2188\\
& 3 & 547.9783 & 3325.1538 & & 2556.8915 & 581.2500 & & 2612.1697 & 602.5000\\
& 4 & 517.0104 & 3141.2211 & & 2417.0151 & 573.4375 & & 2470.1407 & 594.6875\\
& 5 & 496.7651 & 3021.2856 & & 2325.6517 & 567.1875 & & 2377.4376 & 588.4375\\
\cline{2-10}
\multirow{6}{*}{600,300} & 0 & 1188.1381 & 10168.2150 & & 7468.2233 & 648.4375 & & 7524.4775	& 654.6875 \\
& 1 & 359.3716 & 3487.3750 & & 2519.9419 & 496.8750 & & 2556.8817	& 508.1250 \\
& 2 & 262.8568 & 2749.5853 & & 1974.9428 & 466.4063 & &2005.8837	& 477.6553 \\
& 3 & 221.1033 & 2454.6966 && 1755.7921 & 450.7813 & & 1784.3464	& 462.0313 \\
& 4 & 196.0080 & 2286.9560 & & 1630.5261 & 439.8438 & & 1657.7795	 & 451.0938 \\
& 5 & 178.6470 & 2174.9919 & & 1546.5266 & 431.2500 & & 1572.9490 & 442.5000 \\
\hline
\end{tabular}
\label{table:effectoftemperature}
\end{table}

For the rest of the parametric study, the material properties are evaluated at $T=$ 300K for the uniform temperature case. In all the cases, we present the non-dimensionalized critical frequency defined as:
\begin{equation}
\overline{\omega}_{cr}^2 = \omega^2 a^4 \left( \frac{\rho_{mo} h}{D_{mo}} \right)
\end{equation}
where $D_{mo} = \frac{E_mh^3}{12(1-\nu^2)}$ and subscript `o' refers to material properties at $T=$ 300K. Table \ref{table:effectoftemperature} highlights the influence of small aerodynamic damping on the critical aerodynamic pressure and the critical frequency. This study is conduced through the complex eigenvalue analysis. The dynamic pressure corresponding to a particular value of aerodynamic damping $g_\tau = \textup{Img}(\Omega)/\textup{Real(}\Omega)$ is taken as the critical dynamic pressure. It can be seen from Table \ref{table:effectoftemperature} that the aerodynamic damping enhances the value of flutter speed in comparison to the case without aerodynamic damping. The effect of thermal gradient is also given in Table \ref{table:effectoftemperature} by considering appropriate temperature for evaluating the material properties. The temperature is assumed to vary only in the thickness direction and determined by \Eref{eqn:tempsolu}. As expected, the critical pressures and coalescence frequencies decreases under the influence of thermal gradient. The influence of the skewness of the plate on the flutter behaviour is shown in Table \ref{table:effectofskewangle}. It can be seen that with increasing skew angle, the flutter speed increases, while increasing the gradient index decreases the flutter speed. This can be attributed to the resistance offered by the geometry and to the stiffness degradation due to the increase in the metallic volume fraction, respectively. Table \ref{table:effectofboundary} presents the influence of boundary conditions, viz., all edges simply supported and all edges clamped on the flutter characteristics of FGM plate. It can be seen that the critical pressure is more for the clamped plate in comparison with those of simply supported plate as expected. With increasing plate aspect ratio, the critical pressure increases, whilst increasing the material gradient index, the critical pressure decreases. This is true for both the boundary conditions. It can be seen from Tables \ref{table:effectoftemperature}  and Table \ref{table:effectofboundary} that damping and clamped boundary condition can enhance the critical flutter speed. The influence of the plate thickness $a/h$ and the material gradient index $n$ on the flutter characteristics of FGM plates is shown in Table \ref{table:effectofah}. The critical aerodynamic pressure increases with decreasing plate thickness and decreases with material gradient index. Again, the decrease in the critical pressure with increasing material gradient index can be attributed to the stiffness degradation due to increase in the metallic volume fraction.

\begin{table}[htpb]
\renewcommand\arraystretch{1.5}
\caption{Influence of the skewangle $\psi$ on the flutter behaviour of temperature dependent FGM plate (material: Si$_3$N$_4$/SUS304) with $a/b=$ 1, $T_c=$ 600, $a/h=$ 20.}
\centering
\begin{tabular}{lrrrrrr}
\hline
Skew angle & \multicolumn{6}{c}{gradient index, $n$} \\
\cline{2-7}
$\psi$ & 0 & 1 & 2 & 3 & 4 & 5 \\
\hline
0 & 648.4375 & 496.8750 & 466.4063 & 450.7813 & 439.8438 & 431.2500 \\
10$^\circ$&	660.9375 & 507.0313 & 475.0000 & 459.3750 & 448.4375 & 439.0625 \\
15$^\circ$&	678.9063 & 520.3125 & 488.2813	& 471.8750	& 460.1563 & 451.5625 \\
20$^\circ$&	706.2500 & 542.1875 & 508.5938	& 492.1875 & 480.4688 & 471.0938 \\
25$^\circ$&	747.6563 & 573.4375 & 539.0625	& 521.0938 & 509.3750 & 409.2188 \\
30$^\circ$&	806.2500 & 618.7500 & 582.0313	& 563.2813 & 550.7813 & 540.6250 \\
\hline
\end{tabular}
\label{table:effectofskewangle}
\end{table}

\begin{table}[htpb]
\renewcommand\arraystretch{1.5}
\caption{Effect of boundary condition on the flutter behaviour of temperature dependent FGM plate (material: Si$_3$N$_4$/SUS304) with $a/b=$ 1, $T_c=$ 300, $a/h=$ 20 and $\psi=$ 0$^\circ$.}
\centering
\begin{tabular}{lrrrrr}
\hline
gradient index & \multicolumn{2}{c}{SSSS} & & \multicolumn{2}{c}{CCCC} \\
\cline{2-3} \cline{5-6}
$n$ & $a/b=$ 1 & $a/b=$ 2 && $a/b=$ 1& $a/b=$ 2 \\
\hline
0& 476.5625 & 1026.5625	&& 1230.4688	& 2136.7188 \\		
1& 380.4688 & 817.1875	&& 985.9375	& 1713.2813 \\	
2& 364.0625 & 782.0313	&& 942.1875	& 1634.3750 \\	
3& 	357.0313 & 767.1875	&& 922.6563	& 1598.4375 \\
4&	352.3438	 & 57.0313	&&908.5938	 & 1574.2188	\\	
5&	348.4375 & 749.2188	&&898.4375	 & 1555.4688 \\
\hline
\end{tabular}
\label{table:effectofboundary}
\end{table}

\begin{table}[htpb]
\renewcommand\arraystretch{1.5}
\caption{Variation of flutter behaviour with plate aspect ratio $a/h$ for a simply supported FGM square plate with $T_c=$ 300K.}
\centering
\begin{tabular}{lrrrrrr}
\hline
$a/h$ & \multicolumn{6}{c}{gradient index, $n$} \\
\cline{2-7}
& 0 & 1 & 2 & 3 & 4 & 5 \\
\hline
5&	570.3125 & 448.4375 & 427.3438 & 418.7500 & 414.0625 & 410.1563 \\
10&	718.7500	 & 571.8750 & 546.8750 & 535.9375 & 528.9063 & 523.4375 \\
20&	776.5625 & 619.5313 & 592.9688 & 581.2500 & 573.4375 & 567.1875 \\
50&	796.0938 & 635.1563 & 607.8125 & 596.0938 & 587.5000 & 581.2500 \\
100&	798.4375 & 637.5000 & 610.1563 & 598.4375 & 589.8438 & 583.5938 \\
\hline
\end{tabular}
\label{table:effectofah}
\end{table}

\begin{figure}[htpb]
\centering
\includegraphics[scale=0.7]{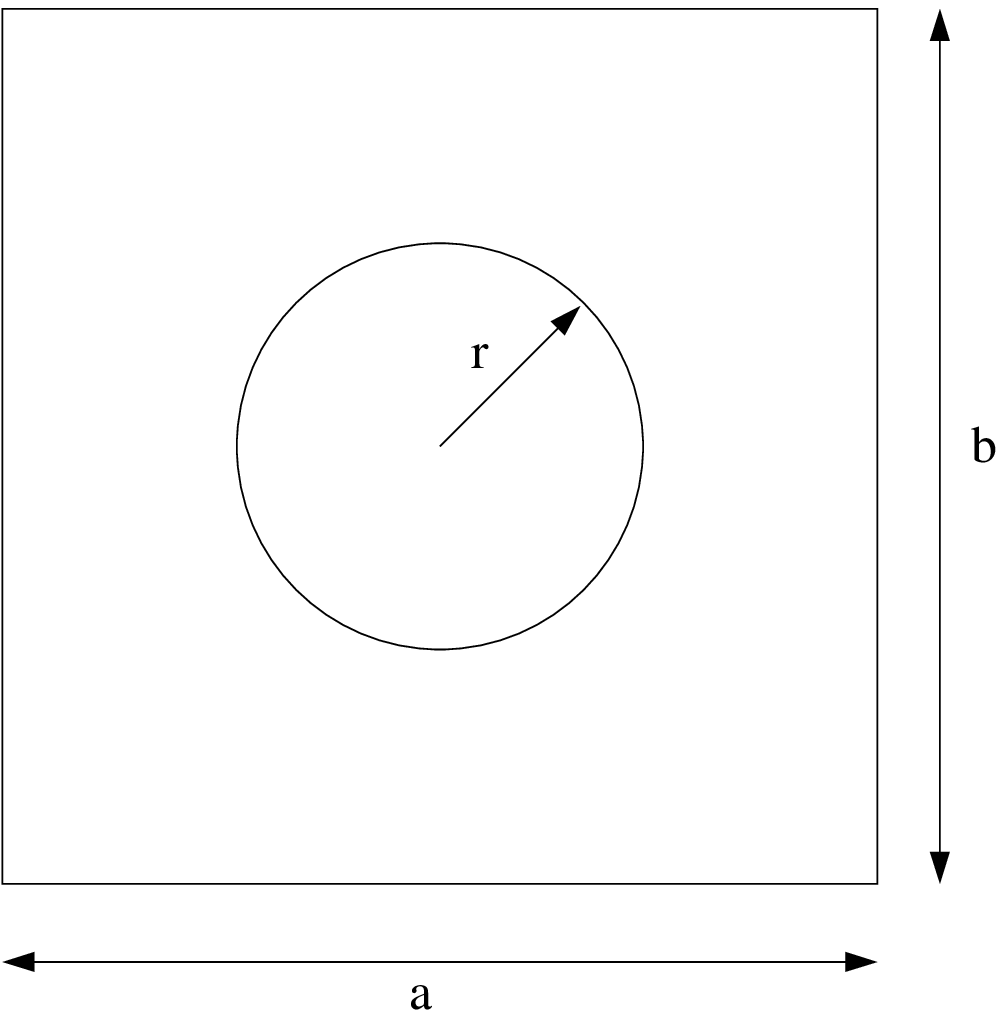}
\caption{Plate with a centrally located circular cutout. $r$ is the radius of the circular cutout.}
\label{fig:phole}
\end{figure}

As a last example, we study the influence of a centrally located circular cutout on the flutter characteristics of FGM plates. A simply supported boundary condition is assumed for this study. \fref{fig:phole} shows the geometry of the plate with a centrally located circular cutout. Table \ref{table:effectofcutout} presents the influence of the size of the centrally located circular cutout on the flutter characteristics of FGM square plate with $a/h=$ 20 and uniform temperature distribution. It can be inferred that increasing the gradient index decreases the critical flutter speed, whilst, increasing the cutout radius, increases the critical flutter speed. This can be attributed to the stiffness degradation due to increase in the metallic volume fraction and due to the presence of a cutout, respectively. It can be seen from the above numerical study that aerodynamic damping, clamped boundary condition and presence of a centrally located circular cutout enhances the critical flutter speed.

\begin{table}[htpb]
\renewcommand\arraystretch{1.5}
\caption{Influence of the radius $r/a$ of a centrally located cutout on the flutter behaviour of temperature dependent simply supported FGM plate (material: Si$_3$N$_4$/SUS304) with $a/b=$ 1, $T_c=$ 300, $a/h=$ 20, $\psi=$ 0$^\circ$.}
\centering
\begin{tabular}{lrrrrrr}
\hline
$r/a$ & \multicolumn{6}{c}{gradient index, $n$} \\
\cline{2-7}
& 0 & 1 & 2 & 3 & 4 & 5 \\
\hline
0	&776.5625	&619.5313	&592.96875&	581.2500&	573.4375&	567.1875\\
0.1	&920.1563 & 	746.5625 	&720.3125	& 707.03125	& 702.65625	& 695.6250\\
0.2	&1727.9688 	&1429.6875	 &1382.6563 &	1361.8750 &	1348.90625	& 1334.8438\\
0.3	& 3179.6875	&2588.2813 &	2464.8438 &	2421.0938 &	2360.1563 &	2338.2813\\
0.4	& 8430.4688	&6798.9063 &	6522.5000&	 6401.2500&	6317.1875 &	6246.8750\\
\hline
\end{tabular}
\label{table:effectofcutout}
\end{table}

\vspace*{0.2cm}
\section{Conclusion}
In this study, a cell based smoothing technique with discrete shear gap method for three-noded triangular element was detailed and used to study the linear flutter characteristics of flat functionally graded material panels. The efficiency and accuracy of the formulation are demonstrated with few numerical examples. From the detailed numerical study, the following can be concluded:
\begin{itemize}
\item With increasing gradient index and the plate aspect ratio, the critical aerodynamic pressure decreases.
\item Damping and clamped boundary condition enhances the critical flutter speed.
\item Increasing the cutout radius, increases the critical aerodynamic pressure when the cutout is centrally located.
\item Thermal gradient decreases the critical flutter speed as expected.
\item Coalescence of higher modes are possible in determining the critical value. This depends on the aspect ratio, the cutout size and the thermal gradient.
\end{itemize}

\section*{Acknowledgements} 
S Natarajan would like to acknowledge the financial support of the School of Civil and Environmental Engineering, The University of New South Wales for his research fellowship since September 2012. 

\section*{References}
\bibliographystyle{plain}
\bibliography{myRefFGM}

\end{document}

%% file: fgm.tex
\subsection{Functionally graded material}
A rectangular plate made of a mixture of ceramic and metal is considered with the coordinates $x,y$ along the in-plane directions and $z$ along the thickness direction (see \fref{fig:platefig}). The material on the top surface $(z=h/2)$ of the plate is ceramic rich and is graded to metal at the bottom surface of the plate $(z=-h/2)$ by a power law distribution. The effective properties of the FGM plate can be computed by using the rule of mixtures or by employing the Mori-Tanaka homogenization scheme. Let $V_i (i=c,m)$ be the volume fraction of the phase material. The subscripts $c$ and $m$ refer to ceramic and metal phases, respectively. The volume fraction of ceramic and metal phases are related by $V_c + V_m = 1$ and $V_c$ is expressed as:
\begin{equation}
V_c(z) = \left( \frac{2z+h}{2h} \right)^n
\end{equation}
where $n$ is the volume fraction exponent $(n \geq 0)$, also known as the gradient index. The variation of the composition of ceramic and metal is linear for $n=$1, the value of $n=$ 0 represents a fully ceramic plate and any other value of $n$ yields a composite material with a smooth transition from ceramic to metal.

\noindent \paragraph*{Mori-Tanaka homogenization method}
Based on the Mori-Tanaka homogenization method, the effective Young's modulus and Poisson's ratio are computed from the effective bulk modulus $K$ and the effective shear modulus $G$ as~\cite{Sundararajan2005}

\begin{align}
\frac{K_{\rm eff}-K_m}{K_c-K_m} =& \frac{V_c}{1 + V_m \frac{3(K_c-K_m)}{3K_m+4G_m}} \nonumber \\
\frac{G_{\rm eff}-G_m}{G_c-G_m} =& \frac{V_c}{1 + V_m \frac{(G_c-G_m)}{(G_m+f_1)} }
\end{align}
where
\begin{equation}
f_1 = \frac{G_m (9K_m+8G_m)}{6(K_m+2G_m)}
\end{equation}
The effective Young's modulus $E_{\rm eff}$ and Poisson's ratio $\nu_{\rm eff}$ can be computed from the following relations:

\begin{equation}
E_{\rm eff} = \frac{9 K_{\rm eff} G_{\rm eff}}{3K_{\rm eff} + G_{\rm eff}}, \hspace{1cm} \nu_{\rm eff} = \frac{3K_{\rm eff} - 2G_{\rm eff}}{2(3K_{\rm eff} + G_{\rm eff})}
\label{eqn:young}
\end{equation}
The effective mass density $\rho$ is computed using the rule of mixtures as $\rho = \rho_c V_c + \rho_m V_m$. The effective heat conductivity $\kappa_{\rm eff}$ and the coefficient of thermal expansion $\alpha_{\rm eff}$ is given by:

\begin{align}
\frac{\kappa_{\rm eff} - \kappa_m}{\kappa_c - \kappa_m} &= \frac{V_c}{1 + V_m \frac{(\kappa_c - \kappa_m)}{3\kappa_m}} \nonumber \\
\frac{\alpha_{\rm eff} - \alpha_m}{\alpha_c - \alpha_m} &= \frac{ \left( \frac{1}{K_{\rm eff}} - \frac{1}{K_m} \right)}{\left(\frac{1}{K_c} - \frac{1}{K_m} \right)}
\label{eqn:thermalcondalpha}
\end{align}

\noindent \paragraph*{Temperature dependent material property} The material properties that are temperature dependent are written as~\cite{Sundararajan2005}:
\begin{equation}
P = P_o(P_{-1}T^{-1} + 1 + P_1T + P_2T^2 + P_3T^3)
\end{equation}
where $P_o,P_{-1},P_1,P_2$ and $P_3$ are the coefficients of temperature $T$ and are unique to each constituent material phase.

\noindent \paragraph*{Temperature distribution through the thickness}
The temperature variation is assumed to occur in the thickness direction only and the temperature field is considered to be constant in the $xy$-plane. In such a case, the temperature distribution along the thickness can be obtained by solving a steady state heat transfer problem:

\begin{equation}
-{d \over dz} \left[ \kappa(z) {dT \over dz} \right] = 0, \hspace{0.5cm} T = T_c ~\textup{at}~ z = h/2;~~ T = T_m ~\textup{at} ~z = -h/2
\label{eqn:heat}
\end{equation}
The solution of \Eref{eqn:heat} is obtained by means of a polynomial series~\cite{wu2004} as
\begin{equation}
T(z) = T_m + (T_c - T_m) \eta(z,h)
\label{eqn:tempsolu}
\end{equation}
where,
\begin{equation}
\begin{split}
\eta(z,h) = {1 \over C} \left[ \left( {2z + h \over 2h} \right) - {\kappa_{cm} \over (n+1)\kappa_m} \left({2z + h \over 2h} \right)^{n+1} + {\kappa_{cm} ^2 \over (2n+1)\kappa_m ^2 } \left({2z + h \over 2h} \right)^{2n+1}
-{\kappa_{cm} ^3 \over (3n+1)\kappa_m ^3 } \left({2z + h \over 2h} \right)^{3n+1} \right. \\ + 
\left. {\kappa_{cm} ^4 \over (4n+1)\kappa_m^4 } \left({2z + h \over 2h} \right)^{4n+1} 
- {\kappa_{cm} ^5 \over (5n+1)\kappa_m ^5 } \left({2z + h \over 2h} \right)^{5n+1} \right] ;
\end{split}
\label{eqn:heatconducres}
\end{equation}

\begin{equation}
\begin{split}
C = 1 - {\kappa_{cm} \over (n+1)\kappa_m} + {\kappa_{cm} ^2 \over (2 n+1)\kappa_m ^2} 
- {\kappa_{cm} ^3 \over (3n+1)\kappa_m ^3} \\ + {\kappa_{cm} ^4 \over (4n+1)\kappa_m ^4}
- {\kappa_{cm} ^5\over (5n+1)\kappa_m ^5}
\end{split}
\end{equation}

%% file: plateTheory.tex
\subsection{Reissner-Mindlin Plates}
The Reissner-Mindlin plate theory, also known as the first order shear deformation theory, takes into account the shear deformation through the thickness, in which the normal to the medium surface remains straight but not necessarily perpendicular to the medium surface. The displacements $u,v,w$ at a point $(x,y,z)$ in the plate (see \fref{fig:platefig}) from the medium surface are expressed as functions of the mid-plane displacements $u_o,v_o,w_o$ and independent rotations $\theta_x,\theta_y$ of the normal in $yz$ and $xz$ planes, respectively, as:

\begin{eqnarray}
u(x,y,z,t) &=& u_o(x,y,t) + z \theta_x(x,y,t) \nonumber \\
v(x,y,z,t) &=& v_o(x,y,t) + z \theta_y(x,y,t) \nonumber \\
w(x,y,z,t) &=& w_o(x,y,t) 
\label{eqn:displacements}
\end{eqnarray}

\begin{figure}[htpb]
\centering
\subfigure[]{\input{./Figures/plate.pstex_t}}
\subfigure[]{\input{./Figures/skew.pstex_t}}
\caption{(a) coordinate system of a rectangular FGM plate, (b) Coordinate system of a skew plate}
\label{fig:platefig}
\end{figure}
where $t$ is the time. The strains in terms of mid-plane deformation can be written as:
\begin{equation}
\bveps  = \left\{ \begin{array}{c} \bveps_p \\ 0 \end{array} \right \}  + \left\{ \begin{array}{c} z \bveps_b \\ \bveps_s \end{array} \right\} 
\label{eqn:strain1}
\end{equation}
The midplane strains $\bveps_p$, the bending strains $\bveps_b$ and the shear strain $\varepsilon_s$ in \Eref{eqn:strain1} are written as:

\begin{eqnarray}
\renewcommand{\arraystretch}{1.2}
\bveps_p = \left\{ \begin{array}{c} u_{o,x} \\ v_{o,y} \\ u_{o,y}+v_{o,x} \end{array} \right\}, \hspace{1cm}
\renewcommand{\arraystretch}{1.2}
\bveps_b = \left\{ \begin{array}{c} \theta_{x,x} \\ \theta_{y,y} \\ \theta_{x,y}+\theta_{y,x} \end{array} \right\}, \nonumber \\
\renewcommand{\arraystretch}{1.2}
\bveps_s = \left\{ \begin{array}{c} \theta _x + w_{o,x} \\ \theta _y + w_{o,y} \end{array} \right\}. \hspace{1cm}
\renewcommand{\arraystretch}{1.2}
\end{eqnarray}
where the subscript `comma' represents the partial derivative with respect to the spatial coordinate succeeding it. The membrane stress resultants $\bn$ and the bending stress resultants $\bfm$ can be related to the membrane strains, $\bveps_p$ and bending strains $\bveps_b$ through the following constitutive relations:

\begin{eqnarray}
\bn &=& \left\{ \begin{array}{c} N_{xx} \\ N_{yy} \\ N_{xy} \end{array} \right\} = \mathbf{A} \bveps_p + \BB \bveps_b \nonumber \\
\bfm &=& \left\{ \begin{array}{c} M_{xx} \\ M_{yy} \\ M_{xy} \end{array} \right\} = \BB \bveps_p + \DD \bveps_b 
\end{eqnarray}
where the matrices $\mathbf{A} = A_{ij}, \BB= B_{ij}$ and $\DD = D_{ij}; (i,j=1,2,6)$ are the extensional, the bending-extensional coupling and the bending stiffness coefficients and are defined as:

\begin{equation}
\left\{ A_{ij}, ~B_{ij}, ~ D_{ij} \right\} = \int_{-h/2}^{h/2} \overline{Q}_{ij} \left\{1,~z,~z^2 \right\}~dz
\end{equation}
Similarly, the transverse shear force $Q = \{Q_{xz},Q_{yz}\}$ is related to the transverse shear strains $\varepsilon_s$ through the following equation:

\begin{equation}
Q_{ij} = E_{ij} \varepsilon_s
\end{equation}
where $E_{ij} = \int_{-h/2}^{h/2} \overline{Q}_{ij} \upsilon_i \upsilon_j~dz;~ (i,j=4,5)$ is the transverse shear stiffness coefficient, $\upsilon_i, \upsilon_j$ are the transverse shear coefficients for non-uniform shear strain distribution through the plate thickness. The stiffness coefficients $\overline{Q}_{ij}$ are defined as:

\begin{eqnarray}
\overline{Q}_{11} = \overline{Q}_{22} = {E(z) \over 1-\nu^2}; \hspace{1cm} \overline{Q}_{12} = {\nu E(z) \over 1-\nu^2}; \hspace{1cm} \overline{Q}_{16} = \overline{Q}_{26} = 0 \nonumber \\
\overline{Q}_{44} = \overline{Q}_{55} = \overline{Q}_{66} = {E(z) \over 2(1+\nu) }
\end{eqnarray}
where the modulus of elasticity $E(z)$ and Poisson's ratio $\nu$ are given by \Eref{eqn:young}. The thermal stress resultant $\bn^{\rm th}$ and the moment resultant $\bfm^{\rm th}$ are:

\begin{eqnarray}
\bn^{\rm th}&=& \left\{ \begin{array}{c} N^{\rm th}_{xx} \\ N^{\rm th}_{yy} \\ N^{\rm th}_{xy} \end{array} \right\} = \int\limits_{-h/2}^{h/2} \overline{Q}_{ij} \alpha(z,T) \left\{ \begin{array}{c} 1 \\ 1 \\ 0 \end{array} \right\} \Delta T(z)~ \rmd z \nonumber \\
\bfm^{\rm th} &=& \left\{ \begin{array}{c} M^{\rm th}_{xx} \\ M^{\rm th}_{yy} \\ M^{\rm th}_{xy} \end{array} \right\} = \int\limits_{-h/2}^{h/2} \overline{Q}_{ij} \alpha(z,T) \left\{ \begin{array}{c} 1 \\ 1 \\ 0 \end{array} \right\} \Delta T(z)~ z ~\rmd z \nonumber \\ 
\end{eqnarray}
where the thermal coefficient of expansion $\alpha(z,T)$ is given by \Eref{eqn:thermalcondalpha} and $\Delta T(z) = T(z)-T_o$ is the temperature rise from the reference temperature and $T_o$ is the temperature at which there are no thermal strains. The strain energy function $U$ is given by:

\begin{equation}
\begin{split}
U(\boldsymbol{\delta}) = {1 \over 2} \int_{\Omega} \left\{ \bveps_p^{\textup{T}} \mathbf{A} \bveps_p + \bveps_p^{\textup{T}} \mathbf{B} \bveps_b + 
\bveps_b^{\textup{T}} \mathbf{B} \bveps_p + \bveps_b^{\textup{T}} \mathbf{D} \bveps_b +  \bveps_s^{\textup{T}} \mathbf{E} \bveps_s - \bveps_b^{\rm T} \bn^{\rm th} - \bveps_b^{\rm T} \bfm^{\rm th} \right\}~ \mathrm{d} \Omega
\end{split}
\label{eqn:potential}
\end{equation}
where $\boldsymbol{\delta} = \{u,v,w,\theta_x,\theta_y\}$ is the vector of the degree of freedom associated to the displacement field in a finite element discretization. Following the procedure given in~\cite{Rajasekaran1973}, the strain energy function $U$ given in~\Eref{eqn:potential} can be rewritten as:

\begin{equation}
U(\boldsymbol{\delta}) = {1 \over 2}  \boldsymbol{\delta}^{\textup{T}} \mathbf{K}  \boldsymbol{\delta}
\label{eqn:poten}
\end{equation}
where $\mathbf{K}$ is the linear stiffness matrix. The kinetic energy of the plate is given by:

\begin{equation}
T(\boldsymbol{\delta}) = {1 \over 2} \int_{\Omega} \left\{p (\dot{u}_o^2 + \dot{v}_o^2 + \dot{w}_o^2) + I(\dot{\theta}_x^2 + \dot{\theta}_y^2) \right\}~\mathrm{d} \Omega
\label{eqn:kinetic}
\end{equation}
where $p = \int_{-h/2}^{h/2} \rho(z)~dz, ~ I = \int_{-h/2}^{h/2} z^2 \rho(z)~dz$ and $\rho(z)$ is the mass density that varies through the thickness of the plate. When the plate is subjected to a temperature field, this in turn results in in-plane stress resultants, $\bn^{\rm th}$. The external work due to the in-plane stress resultants developed in the plate under a thermal load is given by:
\begin{equation}
\begin{split}
V(\boldsymbol{\delta}) = \int\limits_\Omega \left\{ \frac{1}{2} \left[ N_{xx}^{\rm th} w_{,x}^2 + N_{yy}^{\rm th} w_{,y}^2 + 2 N_{xy}^{\rm th}w_{,x}w_{,y}\right] + \right. \\ \left.
\frac{h^2}{24} \left[ N_{xx}^{\rm th} \left( \theta_{x,x}^2 + \theta_{y,x}^2 \right) + N_{yy}^2 \left( \theta_{x,y}^2 + \theta_{y,y}^2 \right) + 2 N_{xy}^{\rm th} \left( \theta_{x,x}\theta_{x,y} + \theta_{y,x}\theta_{y,y} \right) \right] \right\}~ d\Omega
\end{split}
\end{equation}
The work done by the applied non-conservative loads is:
\begin{equation}
W(\boldsymbol{\delta}) = \int_{\Omega} \Delta p w ~\rmd \Omega
\label{eqn:aerowork}
\end{equation}
where $\Delta p$ is the aerodynamic pressure. The aerodynamic pressure based on first-order, high Mach number approximation to linear potential flow is given by:
\begin{equation}
\Delta p = \frac{\rho_a U_a^2}{\sqrt{M_\infty^2 - 1}} \left[ \frac{\partial w}{\partial x} \cos \theta^\prime + \frac{\partial w}{\partial y} \sin \theta^\prime + \left( \frac{1}{U_a} \right) \frac{M_\infty^2 - 2}{M_\infty^2 - 1} \frac{\partial w}{\partial t} \right]
\label{eqn:aeropressure}
\end{equation}
where $\rho_a, U_a, M_\infty$ and $\theta^\prime$ are the free stream air density, velocity of air, Mach number and flow angle, respectively. Substituting \Eref{eqn:poten} - (\ref{eqn:aerowork}) in Lagrange's equations of motion, the following governing equation is obtained:
\begin{equation}
\bfm \ddot{\boldsymbol{\delta}} + g_\tau \mathbf{D}_A \dot{\boldsymbol{\delta}} +  (\KK + \KK_G + \lambda \overline{\mathbf{A}}) \boldsymbol{\delta} = \mathbf{0}
\label{eqn:govereqn}
\end{equation}
where $\KK$ is the stiffness matrix, $\KK_G$ is the geometric stiffness matrix essentially a function of the in-plane stress distribution due to the applied temperature distribution over the plate, $\bfm$ is the consistent mass matrix, $\lambda = \frac{\rho_a U_a^2}{\sqrt{M_\infty^2 - 1}}$, $\overline{\mathbf{A}}$ is the aerodynamic force matrix and $g_\tau = \frac{\lambda (M_\infty^2-1)}{U_a(M^2_\infty-1)}$ is the aerodynamic damping parameter. The damping matrix $\mathbf{D}_A$ can be considered as the scalar multiple of mass matrix by neglecting the shear and rotarty inertia terms of the mass matrix $\bfm$ and after substituting the characteristic of the time function $\ddot{\boldsymbol{\delta}} = -\omega^2 \boldsymbol{\delta}$, the following algebraic equation is obtained:

\begin{equation}
\left[ \left( \KK + \KK_G + \lambda \overline{\mathbf{A}}\right) - \overline{\kappa} \bfm\right] \boldsymbol{\delta} = \mathbf{0}
\label{eqn:finaldiscre}
\end{equation}
where the eigenvalue $\overline{\kappa} = -\omega^2 - g_\tau \omega/(\rho h)$ includes the contribution of aerodynamic damping. \Eref{eqn:finaldiscre} is solved for eigenvalues for a given value of $\lambda$. In the absence of aerodynamic damping, when $\lambda = $0, the eigenvalue of $\omega$ is real and positive, since the stiffness matrix and mass matrix are symmetric and positive definite. However, the aerodynamic matrix $\overline{\mathbf{A}}$ is unsymmetric and hence complex eigenvalues $\omega$ are expected for $\lambda >$ 0. As $\lambda$ increases monotonically from zero, two of these eigenvalues will approach each other and become complex conjugates. In this study, $\lambda_{cr}$ is considered to be the value of $\lambda$ at which the first coalescence occurs. In the presence of aerodynamic damping, the eigenvalues $\overline{\kappa}$, in \Eref{eqn:finaldiscre} becomes complex with increase in the value of $\lambda$. The corresponding frequency can be written as:
\begin{equation}
\overline{\kappa} = -\omega^2 - g_\tau \omega/(\rho h) = \overline{\kappa}_R - i \overline{\kappa}_I
\end{equation}
where the subscripts $R$ and $I$ refer to the real and the imaginary part of the eigenvalue. The flutter boundary is reached $(\lambda = \lambda_{cr})$ when the frequency $\omega$ becomes pure imaginary number, i.e., $\omega = i \sqrt{\overline{\kappa}_R}$ at $g_\tau = \overline{\kappa}_I/\sqrt{\overline{\kappa}_R}$. In practice, the value of $\lambda_{cr}$ is determined from a plot of $\omega_R$ vs $\lambda$ corresponding to $\omega_R = $ 0.




%% file: Figures/plate.pstex_t
\begin{picture}(0,0)%
\includegraphics{./Figures/plate.pstex}%
\end{picture}%
\setlength{\unitlength}{2279sp}%
\begingroup\makeatletter\ifx\SetFigFont\undefined%
\gdef\SetFigFont#1#2#3#4#5{%
  \reset@font\fontsize{#1}{#2pt}%
  \fontfamily{#3}\fontseries{#4}\fontshape{#5}%
  \selectfont}%
\fi\endgroup%
\begin{picture}(7002,4165)(2689,-5264)
\put(4501,-5191){\makebox(0,0)[lb]{\smash{{\SetFigFont{8}{9.6}{\familydefault}{\mddefault}{\updefault}$a$}}}}
\put(8371,-3481){\makebox(0,0)[lb]{\smash{{\SetFigFont{8}{9.6}{\familydefault}{\mddefault}{\updefault}$b$}}}}
\put(9676,-1861){\makebox(0,0)[lb]{\smash{{\SetFigFont{8}{9.6}{\familydefault}{\mddefault}{\updefault}$h$}}}}
\end{picture}%

%% file: Figures/skew.pstex_t
\begin{picture}(0,0)%
\includegraphics{./Figures/skew.pstex}%
\end{picture}%
\setlength{\unitlength}{2279sp}%
\begingroup\makeatletter\ifx\SetFigFont\undefined%
\gdef\SetFigFont#1#2#3#4#5{%
  \reset@font\fontsize{#1}{#2pt}%
  \fontfamily{#3}\fontseries{#4}\fontshape{#5}%
  \selectfont}%
\fi\endgroup%
\begin{picture}(6804,3835)(2659,-5039)
\put(2746,-4156){\makebox(0,0)[lb]{\smash{{\SetFigFont{8}{9.6}{\familydefault}{\mddefault}{\updefault}$\psi$}}}}
\put(2746,-2671){\makebox(0,0)[lb]{\smash{{\SetFigFont{8}{9.6}{\familydefault}{\mddefault}{\updefault}$y$}}}}
\put(5491,-1411){\makebox(0,0)[lb]{\smash{{\SetFigFont{8}{9.6}{\familydefault}{\mddefault}{\updefault}$y^\prime$}}}}
\put(7606,-4471){\makebox(0,0)[lb]{\smash{{\SetFigFont{8}{9.6}{\familydefault}{\mddefault}{\updefault}$x,x^\prime$}}}}
\put(4636,-4966){\makebox(0,0)[lb]{\smash{{\SetFigFont{8}{9.6}{\familydefault}{\mddefault}{\updefault}$a$}}}}
\put(8371,-3301){\makebox(0,0)[lb]{\smash{{\SetFigFont{8}{9.6}{\familydefault}{\mddefault}{\updefault}$b$}}}}
\end{picture}%

%% file: Figures/triangle.pstex_t
\begin{picture}(0,0)%
\includegraphics{./Figures/triangle.pstex}%
\end{picture}%
\setlength{\unitlength}{3947sp}%
\begingroup\makeatletter\ifx\SetFigFont\undefined%
\gdef\SetFigFont#1#2#3#4#5{%
  \reset@font\fontsize{#1}{#2pt}%
  \fontfamily{#3}\fontseries{#4}\fontshape{#5}%
  \selectfont}%
\fi\endgroup%
\begin{picture}(5253,4245)(2161,-4876)
\put(4426,-3361){\makebox(0,0)[lb]{\smash{{\SetFigFont{14}{16.8}{\familydefault}{\mddefault}{\updefault}{\color[rgb]{0,0,0}O}%
}}}}
\put(4201,-811){\makebox(0,0)[lb]{\smash{{\SetFigFont{14}{16.8}{\familydefault}{\mddefault}{\updefault}{\color[rgb]{0,0,0}1}%
}}}}
\put(2176,-4111){\makebox(0,0)[lb]{\smash{{\SetFigFont{14}{16.8}{\familydefault}{\mddefault}{\updefault}{\color[rgb]{0,0,0}2}%
}}}}
\put(7276,-4861){\makebox(0,0)[lb]{\smash{{\SetFigFont{14}{16.8}{\familydefault}{\mddefault}{\updefault}{\color[rgb]{0,0,0}3}%
}}}}
\put(3601,-2761){\makebox(0,0)[lb]{\smash{{\SetFigFont{14}{16.8}{\familydefault}{\mddefault}{\updefault}{\color[rgb]{0,0,0}$\Delta_1$}%
}}}}
\put(3976,-3886){\makebox(0,0)[lb]{\smash{{\SetFigFont{14}{16.8}{\familydefault}{\mddefault}{\updefault}{\color[rgb]{0,0,0}$\Delta_2$}%
}}}}
\put(5026,-3061){\makebox(0,0)[lb]{\smash{{\SetFigFont{14}{16.8}{\familydefault}{\mddefault}{\updefault}{\color[rgb]{0,0,0}$\Delta_3$}%
}}}}
\end{picture}%

%% file: Figures/stdtri.pstex_t
\begin{picture}(0,0)%
\includegraphics{./Figures/stdtri.pstex}%
\end{picture}%
\setlength{\unitlength}{3947sp}%
\begingroup\makeatletter\ifx\SetFigFont\undefined%
\gdef\SetFigFont#1#2#3#4#5{%
  \reset@font\fontsize{#1}{#2pt}%
  \fontfamily{#3}\fontseries{#4}\fontshape{#5}%
  \selectfont}%
\fi\endgroup%
\begin{picture}(6102,4887)(2911,-5236)
\put(4816,-4021){\makebox(0,0)[lb]{\smash{{\SetFigFont{14}{16.8}{\familydefault}{\mddefault}{\updefault}{\color[rgb]{0,0,0}$\xi$}%
}}}}
\put(3451,-4786){\makebox(0,0)[lb]{\smash{{\SetFigFont{14}{16.8}{\familydefault}{\mddefault}{\updefault}{\color[rgb]{0,0,0}1}%
}}}}
\put(7876,-3286){\makebox(0,0)[lb]{\smash{{\SetFigFont{14}{16.8}{\familydefault}{\mddefault}{\updefault}{\color[rgb]{0,0,0}2}%
}}}}
\put(4726,-1411){\makebox(0,0)[lb]{\smash{{\SetFigFont{14}{16.8}{\familydefault}{\mddefault}{\updefault}{\color[rgb]{0,0,0}3}%
}}}}
\put(4126,-1471){\makebox(0,0)[lb]{\smash{{\SetFigFont{14}{16.8}{\familydefault}{\mddefault}{\updefault}{\color[rgb]{0,0,0}d}%
}}}}
\put(2926,-3001){\makebox(0,0)[lb]{\smash{{\SetFigFont{14}{16.8}{\familydefault}{\mddefault}{\updefault}{\color[rgb]{0,0,0}c}%
}}}}
\put(5611,-5221){\makebox(0,0)[lb]{\smash{{\SetFigFont{14}{16.8}{\familydefault}{\mddefault}{\updefault}{\color[rgb]{0,0,0}b}%
}}}}
\put(7906,-4036){\makebox(0,0)[lb]{\smash{{\SetFigFont{14}{16.8}{\familydefault}{\mddefault}{\updefault}{\color[rgb]{0,0,0}b}%
}}}}
\put(4021,-3706){\makebox(0,0)[lb]{\smash{{\SetFigFont{14}{16.8}{\familydefault}{\mddefault}{\updefault}{\color[rgb]{0,0,0}$\eta$}%
}}}}
\end{picture}%

%% file: CSFEM_Flut.bbl
\begin{thebibliography}{10}

\bibitem{aliatwal1980}
R~Ali and SJ~Atwal.
\newblock Prediction of natural frequencies of vibration of rectangular plates
  with rectangular cutouts.
\newblock {\em Computers and Structures}, 12:819--823, 1980.

\bibitem{baiznatarajan2011}
Pedro~M Baiz, S~Natarajan, SPA Bordas, P~Kerfriden, and T~Rabczuk.
\newblock Linear buckling analysis of cracked plates by {SFEM and XFEM}.
\newblock {\em Journal of Mechanics of Materials and Structure}, 6:1213--1238,
  2011.

\bibitem{bletzingerbischoff2000}
KU~Bletzinger, M~Bischoff, and E~Ramm.
\newblock A unified approach for shear locking free triangular and rectangular
  shell finite elements.
\newblock {\em Computers and Structures}, 75:321--334, 2000.

\bibitem{bordasnatarajan2010}
S~Bordas and S~Natarajan.
\newblock On the approximation in the smoothed finite element method {(SFEM)}.
\newblock {\em International Journal for Numerical Methods in Engineering},
  81:660--670, 2010.

\bibitem{chowdarysinha1996}
TVR Chowdary, PK~Sinha, and S~Parthan.
\newblock Finite elemen flutter analysis of composite skew panels.
\newblock {\em Computers \& Structures}, 58:613--620, 1996.

\bibitem{Ferreira2006}
AJM Ferreira, RC~Batra, CMC Roque, LF~Qian, and RMN Jorge.
\newblock Natural frequencies of functionally graded plates by a meshless
  method.
\newblock {\em Composite Structures}, 75:593--600, 2006.

\bibitem{ganapathiprakash2006}
M~Ganapathi, T~Prakash, and N~Sundararajan.
\newblock Influence of functionally graded material on buckling of skew plates
  under mechanical loads.
\newblock {\em {ASCE} Journal of Engineering Mechanics}, 132:902--905, 2006.

\bibitem{Huang2011}
CS~Huang, OG~McGee III, and MJ~Chang.
\newblock Vibrations of cracked rectangular {FGM} thick plates.
\newblock {\em Composite Structures}, 93(7):1747--1764, 2011.

\bibitem{huangsakiyama1999}
M~Huang and T~Sakiyama.
\newblock Free vibration analysis of rectangular plates with variously-shaped
  holes.
\newblock {\em Journal of Sound and Vibration}, 226(4):769--786, 1999.

\bibitem{ibrahimtawfik2007}
HH~Ibrahim, M~Tawfik, and M~Al-Ajmi.
\newblock Thermal buckling and nonlinear flutter behavior of functionally
  graded material panels.
\newblock {\em Journal of Aircraft}, 44:1610--1617, 2007.

\bibitem{janghorbanzare2011}
Maziar Janghorban and Amin Zare.
\newblock Thermal effect on free vibration analysis of functionally graded
  arbitrary straight-sided plates with different cutouts.
\newblock {\em Latin American Journal of Solids and Structures}, 8:245--257,
  2011.

\bibitem{jhakant2013}
DK~Jha, Tarun Kant, and RK~Singh.
\newblock A critical review of recent research on functionally graded plates.
\newblock {\em Composite Structures}, 96:833--849, 2013.

\bibitem{Kitipornchai2009}
S~Kitipornchai, LL~Ke, and J~Yang andY Xiang.
\newblock Nonlinear vibration of edge cracked functionally graded {T}imoshenko
  beams.
\newblock {\em Journal of Sound and Vibration}, 324:962--982, 2009.

\bibitem{liaosun1993}
C.-L Lian and YW~Sun.
\newblock Flutter analysis of stiffened laminated composite plates and shells
  in supersonic flow.
\newblock {\em AIAA J}, 31:1897--1905, 1993.

\bibitem{liudai2007}
GR~Liu, KY~Dai, and TT~Nguyen.
\newblock A smoothed finite elemen for mechanics problems.
\newblock {\em Computational Mechanics}, 39:859--877, 2007.

\bibitem{natarajanbaiz2011a}
S~Natarajan, PM~Baiz, SPA Bordas, P~Kerfriden, and T~Rabczuk.
\newblock Natural frequencies of cracked functionally graded material plates by
  the extended finite element method.
\newblock {\em Composite Structures}, 93:3082--3092, 2011.

\bibitem{natarajanbaiz2011}
S~Natarajan, PM~Baiz, M~Ganapathi, P~Kerfriden, and S~Bordas.
\newblock Linear free flexural vibration of cracked functionally graded plates
  in thermal environment.
\newblock {\em Computers and Structures}, 89:1535--1546, 2011.

\bibitem{natarajanmanickam2012}
S~Natarajan and Ganapathi Manickam.
\newblock Bending and vibration of functionally graded material sandwich plates
  using an accurate theory.
\newblock {\em Finite Elements in Analysis and Design}, 57:32--42, 2012.

\bibitem{navazihaddadpour2007}
HM~Navazi and H~Haddadpour.
\newblock Aero-thermoelastic stability of functionally graded plates.
\newblock {\em Composite Structures}, 80:580--586, 2007.

\bibitem{paramasivam1973}
P~Paramasivam.
\newblock Free vibration of square plates with square openings.
\newblock {\em Journal of Sound and Vibration}, 30:173--178, 1973.

\bibitem{prakashganapathi2006}
T~Prakash and M~Ganapathi.
\newblock Supersonic flutter characteristics of functionally graded flat panels
  including thermal effects.
\newblock {\em Composite Structures}, 72:10--18, 2006.

\bibitem{Qian2004a}
L.~C. Qian, R.~C. Batra, and L.~M. Chen.
\newblock Static and dynamic deformations of thick functionally graded elastic
  plates by using higher order shear and normal deformable plate theory and
  meshless local {P}etrov {G}alerkin method.
\newblock {\em Composites Part B: Engineering}, 35:685--697, 2004.

\bibitem{R2010}
Ahmad~Akbari R, Akbar Bagri, St\'ephane Bordas, and Timon Rabczuk.
\newblock Analysis of thermoelastic waves in a two-dimensional functionally
  graded materials domain by the meshless local {P}etrov-{G}alerkin method.
\newblock {\em Computer Modelling in Engineering and Science}, 65:27--74, 2010.

\bibitem{rahimabadinatarajan2013}
AA~Rahimabadi, S~Natarajan, and S~Bordas.
\newblock Vibration of functionally material plates with cutouts and cracks in
  thermal environment.
\newblock {\em Advances in Crack Growth Modeling, Key Engineering Materials},
  560:157--180, 2013.

\bibitem{Rajasekaran1973}
S~Rajasekaran and DW~Murray.
\newblock Incremental finite element matrices.
\newblock {\em {ASCE} Journal of Structural Divison}, 99:2423--2438, 1973.

\bibitem{reddy1982}
JN~Reddy.
\newblock Large amplitude flexural vibration of layered composite plates with
  cutouts.
\newblock {\em Journal of Sound and Vibration}, 83(1):1--10, 1982.

\bibitem{Reddy2000}
JN~Reddy.
\newblock Analysis of functionally graded plates.
\newblock {\em International Journal for Numerical Methods in Engineering},
  47:663--684, 2000.

\bibitem{Reddy1998}
JN~Reddy and CD~Chin.
\newblock Thermomechanical analysis of functionally graded cylinders and
  plates.
\newblock {\em Journal of Thermal Stresses}, 21:593--629, 1998.

\bibitem{Singh2011}
MK~Singh, T~Prakash, and M~Ganapathi.
\newblock Finite element analysis of functionally graded plates under
  transverse load.
\newblock {\em Finite Elements in Analysis and Design}, 47:453--460, 2011.

\bibitem{singhaganapathi2005}
MK~Singha and M~Ganapathi.
\newblock A parametric study on supersonic flutter behavior of laminated
  composite skew panels.
\newblock {\em Composite Structures}, 69:55--63, 2005.

\bibitem{sivakumariyengar1998}
K~Sivakumar, NGR Iyengar, and Kalyanmoy Deb.
\newblock Optimum design of laminated composite plates with cutouts using a
  genetic algorithm.
\newblock {\em Composite Structures}, 42:265--279, 1998.

\bibitem{Sundararajan2005}
N~Sundararajan, T~Prakash, and M~Ganapathi.
\newblock Nonlinear free flexural vibrations of functionally graded rectangular
  and skew plates under thermal environments.
\newblock {\em Finite Elements in Analysis and Design}, 42:152--168, 2005.

\bibitem{valizadehnatarajan2013}
Navid Valizadeh, Sundararajan Natarajan, Octavio~A Gonzalez-Estrada, Timon
  Rabczuk, Tinh~Quoc Bui, and St\'ephane~PA Bordas.
\newblock Nurbs-based finite element analysis of functionally graded plates:
  static bending, vibration, buckling and flutter.
\newblock {\em Composite Structures}, 99:309--326, 2013.

\bibitem{wu2004}
L~Wu.
\newblock Thermal buckling of a simply supported moderately thick rectangular
  {FGM} plate.
\newblock {\em Composite Structures}, 64:211--218, 2004.

\bibitem{Yang2010}
J~Yang, YX~Hao, W~Zhang, and S~Kitipornchai.
\newblock Nonlinear dynamic response of a functionally graded plate with a
  through-width surface crack.
\newblock {\em Nonlinear {D}ynamics}, 59:207--219, 2010.

\bibitem{Yang2002}
J~Yang and H.-S Shen.
\newblock Vibration characteristic and transient response of shear-deformable
  functionally graded plates in thermal environment.
\newblock {\em Journal of Sound and Vibration}, 255:579--602, 2002.

\end{thebibliography}
